\documentclass[11pt,a4paper]{amsart}

\usepackage{amsmath}
\usepackage{amssymb}
\usepackage{latexsym}
\usepackage{xypic}
\usepackage{amscd,amsthm}

\newcommand{\prf}{{\em Proof}. }

\input{epsf}

\newcommand{\str}{^\prime}

\newcommand{\card}[1]{{\mid\! #1 \!\mid}}

\def\R{{\mathbb R}}
\def\Z{{\mathbb Z}}

\def\V{{\mathcal V}}

\def\conv{{\rm conv}}

\def\Vol{{\rm Vol}}

\def\Sec{{\rm Sec}}
\def\Secf{\Sec(\Delta)}
\def\Discr{{\rm Discr}}
\def\Res{{\rm Res}}
\def\deg{{\rm deg}}

\def\MR{{M_{\R}}}

\def\Deltas{\Delta\str}

\def\intr{{\rm int}}

\newtheorem*{theorem*}{Theorem}
\newtheorem*{corollary*}{Corollary}
\newtheorem{theorem}{Theorem}[section]
\newtheorem{definition}[theorem]{Definition}
\newtheorem{remark}[theorem]{Remark}

\newtheorem{corollary}[theorem]{Corollary}
\newtheorem{proposition}[theorem]{Proposition}
\newtheorem{lemma}[theorem]{Lemma}

\newtheorem{question}[theorem]{Question}

\begin{document}
\title[Multiples of lattice polytopes]{Multiples of lattice 
polytopes without interior lattice points}

\author[Victor Batyrev]{Victor Batyrev}
\address{Mathematisches Institut, Universit\"at T\"ubingen, Auf der 
Morgenstelle 10, 72076 T\"ubingen, Germany}
\email{victor.batyrev@uni-tuebingen.de}

\author[Benjamin Nill]{Benjamin Nill}
\address{Arbeitsgruppe Gitterpolytope, Freie Universit\"at Berlin, Arnimallee 3, 14195 Berlin, Germany}
\email{ nill@math.fu-berlin.de}

\begin{abstract}
Let $\Delta$ be an $n$-dimensional lattice polytope. The smallest non-negative 
integer $i$ such that $k\Delta$ contains no interior lattice 
points for $1 \leq k \leq n - i$ we call the degree of $\Delta$. 
We consider  
lattice polytopes of fixed degree $d$ and arbitrary dimension $n$. Our main 
result is a complete classification of $n$-dimensional 
lattice polytopes of degree $d=1$. This is a generalization of the 
classification of lattice polygons $(n=2)$ without interior lattice points 
due to Arkinstall, Khovanskii, Koelman and \linebreak Schicho. Our classification shows that 
the secondary polytope ${\rm Sec}(\Delta)$ of a lattice polytope of 
degree $1$ is always a simple polytope. 
\end{abstract}



\maketitle

\section{Introduction}

Let $M \cong \Z^n$ be an $n$-dimensional lattice and 
$\Delta \subseteq \MR \cong M \otimes_\Z \R$ an $n$-dimensional 
 lattice polytope (i.e., the set of all 
vertices $\V(\Delta)$ of $\Delta$ is contained in $M$).

It is well-known (cf. \cite{Sta80,Sta86}) that the power series 
\[ P(\Delta, t) :=  \sum\limits_{k \geq 0} \card{(k \Delta) \cap M} \, t^k \]
is a rational function of the following form:
\[ P(\Delta,t) =  \frac{h^*_0 + h^*_1 t + \cdots + h^*_n t^n}{(1-t)^{n+1}}, \]
where $h^*_0, \ldots, h^*_n$ are non-negative integers 
satisfying the conditions $h^*_0 =1$, $h^*_1 = \card{\Delta \cap M} - n - 1$, and $$h^*_0 + \cdots + h^*_n 
=\Vol(\Delta) = n! \{ \mbox{\rm volume  of }\Delta\}. $$  

\begin{definition}{\rm The polynomial 
$ \sum_i h^*_i t^i := (1 -t)^{n+1} P(\Delta,t)$ 
we call the {\em $h^*$-polynomial} of an 
$n$-dimensional 
lattice polytope $\Delta$ and denote it by $h^*_\Delta$. 
The degree of $h^*_\Delta$ we call the {\em degree} 
of $\Delta$ and denote it by  
$\deg(\Delta)$.}
\end{definition}

\begin{remark}{\rm 
Let 
\[ Q(\Delta,t):= \sum\limits_{k \geq 0} \card{\intr(k \Delta) \cap M} \, t^k. \]
Then by the Ehrhart reciprocity theorem (cf. \cite{Sta80,Sta86}) 
one has
\[ Q(\Delta,t) = \frac{h^*_n t  + \cdots + h^*_0 t^{n+1}}{(1-t)^{n+1}}.\]
If $d:= \deg(\Delta)$, then  
\[ h^*_n = h^*_{n-1} = \cdots =h^*_{d+1} =0, \;\;
h^*_d \neq 0 \] 
and the power series $Q(\Delta,t)$ can be factored as 
$h^*_d t^{n-d +1} (1 + \sum_{i>0} c_it^i)$, i.e., we obtain
\[  \card{\intr(\Delta) \cap M} =  \card{\intr(2\Delta) \cap M} = 
\cdots  =  \card{\intr((n-d)\Delta) \cap M} =0, \]
\[ 
   h^*_d = \card{\intr((n-d+1)\Delta) \cap M} \neq 0. \]
Thus  the number  $\deg(\Delta)$  can be defined also as a 
smallest non-negative integer $i$ such that $k\Delta$ contains no 
interior lattice 
points for $1 \leq k \leq n - i$.
}
\label{degree2}
\end{remark}

\begin{remark} 
{\rm We define   $\deg(\Delta)$ for an arbitrary  
lattice polytope $\Delta \subset \MR$ of dimension $\leq n$ 
as the degree of the ($h^*$-)polynomial 
\[ (1-t)^{\dim(\Delta) + 1} P(\Delta, t). \]
}
\label{degree3} 
\end{remark} 

Since for any lattice polytope one has $\deg (\Delta) \leq \dim(\Delta)$, 
it is interesting to ask what can be said about $\Delta$ when 
$\deg(\Delta) \ll \dim(\Delta)$. Our observation is that in this case 
$\deg (\Delta)$ can be considered as a ``lattice dimension'' of $\Delta$. 
In particular, we will show that many examples of lattice polytopes of 
small degree $d$ can be constructed  from lattice polytopes of small dimension.

The following statement immediately follows from the definition 
of $\deg(\Delta)$: 

\begin{proposition}
For a lattice polytope $\Delta$ the following statements are equivalent:
\begin{enumerate}
\item $\deg(\Delta) = 0$;
\item $\Vol(\Delta) = 1$;
\item $\Delta$ is a {\em basic} simplex, i.e.,
 the vertices form an affine lattice basis.
\end{enumerate}
\label{basis}
\end{proposition}

Our main purpose will be a complete classification of 
lattice polytopes of degree $1$ (see Section 2). 
We immediately get from the definition of $\deg(\Delta)$ and the equation $h^*_1 = \card{\Delta \cap M} - n - 1$: 

\begin{proposition}
Let $\Delta$ be an $n$-dimensional lattice polytope. 
Then the following statements are equivalent:

\begin{enumerate}
\item $\deg(\Delta) \leq 1$;
\item $\card{\Delta \cap M} = \Vol(\Delta) + n$.
\end{enumerate}
\label{ele}
\end{proposition}

First we notice that the condition $\deg(\Delta) \ll \dim(\Delta)$
puts rather strong restrictions on the combinatorics of $\Delta$:

\begin{proposition}
Let $\Delta$ be an $n$-dimensional  lattice polytope of degree $d$. 
Then any $k$ lattice points in $\Delta$ such that $k \leq  n -d$ 
are contained in a proper face of $\Delta$. In particular, any two 
vertices of $\Delta$ are contained in a proper face of $\Delta$, if $n-d\geq 2$. 
\end{proposition}

\noindent
\prf Assume that   $v_1, \ldots, v_k$ are  lattice points in $\Delta$ 
which are not contained in a proper face of $\Delta$. Then 
$v_1 + \cdots + v_k \in \intr(k\Delta)$. By \ref{degree2}, we obtain 
$d= \deg(\Delta) \geq n +1 - k$, i.e., $k \geq n-d+1$. Contradiction.

\hfill $\Box$ 

The degree of a lattice polytope is a monotone function, 
this follows directly from the so-called {\em monotonicity theorem} of Stanley:

\begin{theorem}[Stanley \cite{Sta93}] 
Let $\Delta'$ be a lattice subpolytope of a lattice polytope
$\Delta$, i.e., $\Delta'$ is a convex hull of finitely many lattice 
points 
in $\Delta$. If $h^*_{\Delta'} = \sum_i {h'}^*_i t^i$ 
and $h^*_\Delta =  \sum_i h^*_i t^i$, then 
${h'}^*_i \leq h^*_i$ for all $i$.
\end{theorem}

\begin{corollary} 
Let $\Delta'$ be a lattice subpolytope of a lattice polytope $\Delta$. Then 
\[ \deg(\Delta') \leq \deg(\Delta). \]
\label{monot}
\end{corollary}

\begin{corollary}
Let $\Delta$ be an $n$-dimensional  lattice polytope of degree $d$. 
Then every lattice point in $\Delta \cap M$ is contained 
in a face of $\Delta$ of dimension $\leq d$. 
\label{lat-poin}
\end{corollary}

\noindent
\prf Assume that there exists a lattice point $v \in \Delta \cap M$ 
which is contained in the relative interior of 
a $m$-dimensional face $\Theta$  of $\Delta$. Then $\deg(\Theta) 
= \dim(\Theta) =m$. By \ref{monot}, $m = \deg(\Theta) \leq \deg(\Delta) =d$. 
\hfill $\Box$ 
\bigskip 

\begin{proposition}
Let $M \cong  
M' \oplus M''$ be a splitting of an $n$-dimensional lattice into 
direct sum of two sublattices of dimensions $n'$ and $n''$. 
Denote by $\psi$ the canonical surjective homomorphism $\psi\; : \; M \to M''$ 
with the kernel $M'$. Let $\Delta \subset \MR$ be an $n$-dimensional 
lattice polytope and $\psi(\Delta) \subset \MR''$ its $n''$-dimensional 
image in $\MR''$. Then 
\[ \deg(\Delta) \leq \deg(\psi(\Delta)) + n'. \]
\label{proj}
\end{proposition}

\begin{proof} 
Since $\psi( \intr(k\Delta) \cap M) \subset \intr(k\psi(\Delta)) \cap M''$, 
we obtain that 
\[ \card{ \intr(k\Delta) \cap M } \neq 0 \Rightarrow   
\card{ \intr(k\psi(\Delta)) \cap M'' } \neq 0. \]
For $k = n - \deg(\Delta) +1$, by  \ref{degree2}, the latter implies 
$$\deg(\psi(\Delta)) \geq \dim(\psi(\Delta)) - ( n - \deg(\Delta)) 
= \deg(\Delta) - n'.$$ 
\end{proof}

\begin{definition} {\rm Let $M'$ be a $n'$-dimensional lattice. 
Consider $r+1$ 
lattice polytopes $\Delta_0, \ldots, \Delta_r \subset 
\MR' \cong M' \otimes_{\Z} \R$ and set  
$M'':= \Z^{r}$. Denote by $e_0, \ldots e_{r}$ the vertices of the standard 
basic lattice simplex in $\R^r$. The convex hull of  $(\Delta_0, e_0), 
\ldots, (\Delta_r, e_r) $ in $\MR:= \MR' \oplus \MR''$ is called the 
{\em Cayley polytope} of 
$\Delta_0, \ldots, \Delta_r$ and will be denoted by 
$\Delta_0 * \cdots * \Delta_r$. 
One has 
$$\dim (\Delta_0 * \cdots * \Delta_r) = 
r + \dim (  \Delta_0 + \cdots + \Delta_r), $$
If $\dim(\Delta_1) = \cdots = \dim(\Delta_r) =0$, then 
$\Delta_0 * \cdots * \Delta_r$ will be called the {\em $r$-fold pyramid over 
} $\Delta_0$. }
\end{definition}

\begin{proposition}
In the above notation, one has 
\[ \deg(\Delta_0 * \cdots * \Delta_r) \leq 
\dim (\Delta_0 + \cdots + \Delta_r) \leq n'. \]
Moreover, if $\Delta_0 * \cdots * \Delta_r$ is an $r$-fold pyramid over 
$\Delta_0$, then 
\[   \deg(\Delta_0 * \cdots * \Delta_r) =  \deg(\Delta_0). \]
\label{Cay}
\end{proposition}

\noindent
\prf Let $\Delta:=\Delta_0 * \cdots * \Delta_r$. Without loss of generality
we can assume that $n' = \dim (  \Delta_0 + \cdots + \Delta_r)$ and therefore 
$\dim(\Delta) = n' + r = \dim_\R \MR$. 
If $\psi\;:\; M \to M'' = \Z^r$ is the canonical projection, then 
$\psi(\Delta)$ is the standard basic lattice 
simplex in $\R^r$. By \ref{basis}, $\deg(\psi(\Delta)) =0$. Hence, by 
\ref{proj}, 
$\deg(\Delta) \leq n'$. 

If $\Delta$ is an $r$-fold pyramid over $\Delta$, then 
$P(\Delta_0, t) = (1-t)^r P(\Delta, t)$ and $\dim(\Delta) = \dim(\Delta_0) + r$. 
Therefore the polynomials  $(1-t)^{\dim(\Delta_0) + 1} P(\Delta_0, t)$ 
and  $(1-t)^{\dim(\Delta) + 1} P(\Delta, t)$ are the same. By \ref{degree3}, 
$ \deg(\Delta_0) =  \deg(\Delta)$. 

\hfill $\Box$

From \ref{Cay} we see that Cayley polytopes give many examples of 
lattice polytopes of small degree $d$ and arbitrary large dimension $n$. 
In this connection, it would be interesting to know whether the following 
converse version of \ref{Cay} ist true: 

\begin{question} 
Let $d$ be a positive integer. Does there exist a constant $N(d)$ such 
that every lattice polytope $\Delta$ of degree $d$ and dimension $n \geq N(d)$
is a Cayley polytope $\Delta_0 * \Delta_1$ of some lattice polytopes 
$\Delta_0$ and $\Delta_1$ with $\dim(\Delta_0), \dim(\Delta_1) < n$.  
\end{question}     

We show below that for $d =1$ the answer is positive, if we take $N(1) = 3$.

\section{Main theorem on lattice polytopes of degree 1}

The following definition is inspired by the notion 
{\em Lawrence polytope} (see \cite{HRS00}):

\begin{definition}{\rm 
We call an $n$-dimensional lattice polytope $\Delta$ $(n \geq 1)$ 
a {\em Lawrence prism} with {\em heights} $h_1, \ldots, h_n$, 
if there exists an affine lattice basis $e_0, \ldots, e_n$ of $M$ 
and non-negative integers $h_1, \ldots, h_n$ 
such that 
$$\Delta = 
\conv(e_0,e_0 + h_1 (e_n - e_0), e_1, e_1 + h_2 (e_n-e_0), \ldots, e_{n-1}, 
e_{n-1} + h_n (e_n-e_0)).$$ 
In this case, the vector $e_n - e_0$ is called a {\em direction} of $\Delta$.
A Lawrence prism can be considered as  the {\em Cayley polytope} of 
$n$ segments $$[0,h_1], \ldots, [0,h_n] \subset \R.$$
}
\end{definition}

\begin{definition}{\rm
We call an $n$-dimensional lattice polytope  $\Delta$ $(n \geq 2)$ 
{\em exceptional}, 
if there exists an affine lattice basis $e_0, \ldots, e_n$ of $M$ 
such that $$\Delta = 
\conv(e_0,e_0 + 2 (e_1 - e_0), e_0 + 2 (e_2 - e_0), e_3, \ldots, e_n).$$ 
In particular, $\Delta$ is a simplex which is the $(n-2)$-fold  pyramid 
over the 2-dimensional basic simplex multiplied by $2$. 
\label{defi}
}
\end{definition}

The following two figures show  
the 2-dimensional Lawrence prism 
with the heights $h_1 = 3$, $h_2 = 2$, and the exceptional triangle:

\centerline{\epsfxsize1.5in\epsffile{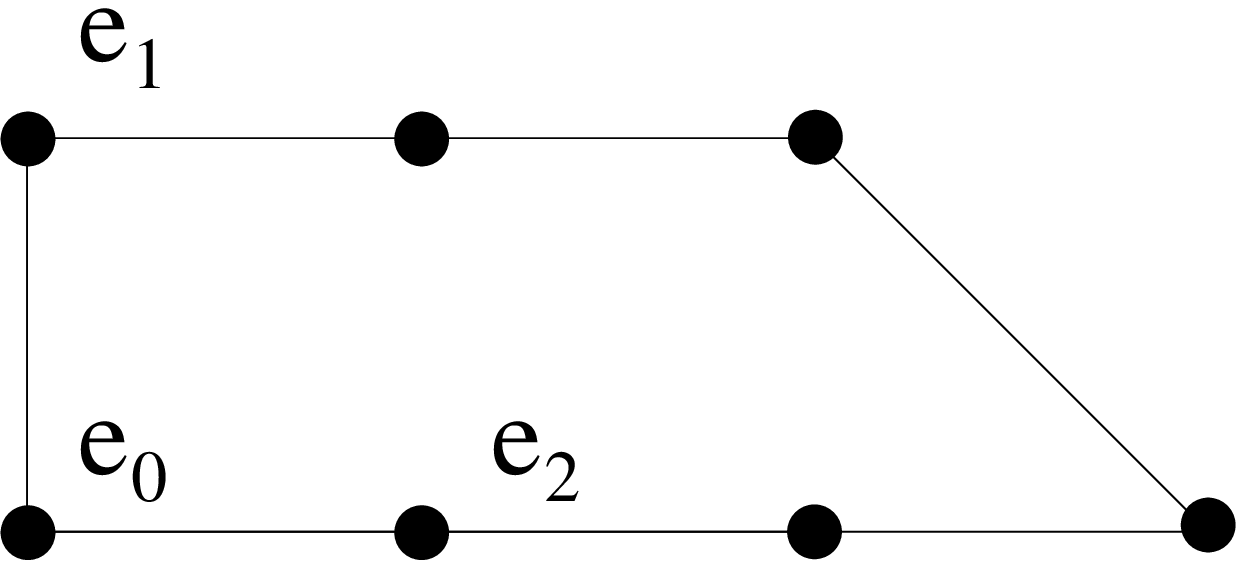}\quad\quad\quad \epsfxsize1in\epsffile{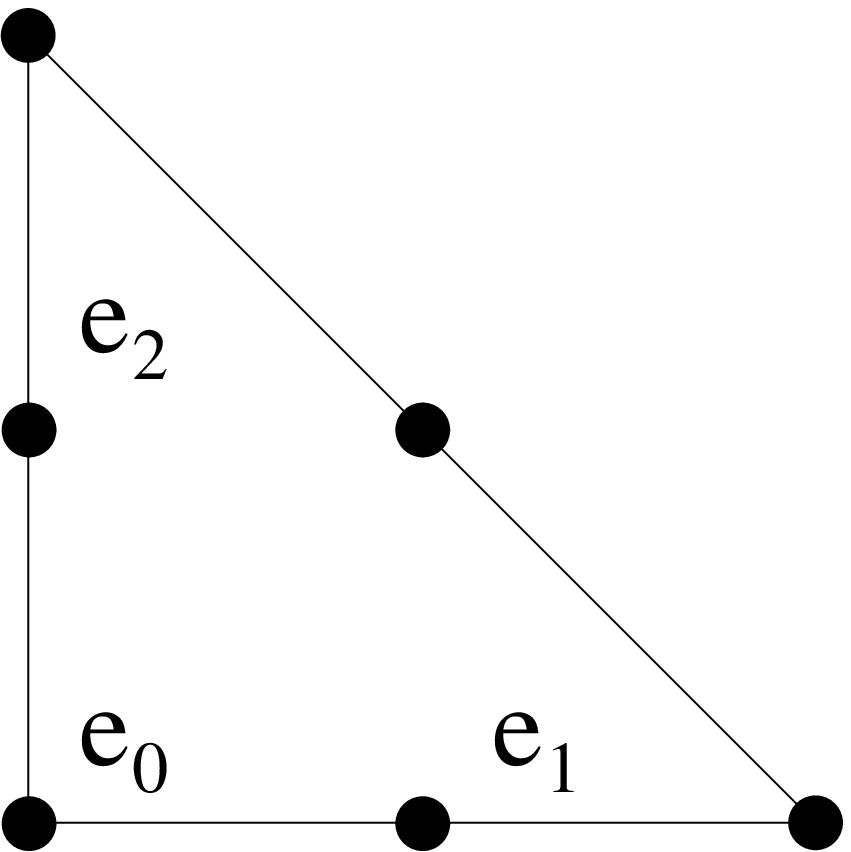}}

\begin{remark} 
{\rm 
\noindent
\begin{enumerate} 
\item We note that an exceptional simplex can never 
be a Lawrence prism, because  it has  two intersecting edges 
each containing more than two lattice points. 

\item For any fixed dimension $n$ there exists exactly one exceptional 
simplex up to isomorphism. On the other hand, $n$-dimensional  
Lawrence prisms form a countably infinite family for any fixed $n$. 

\item We can consider the basic 
lattice simplex of degree $0$ as a particular case 
of a Laurence polytope with the 
heights $h_1 =1$ and $h_2 = \cdots = h_n =0$.   
\end{enumerate}}
\end{remark}

\begin{proposition}
An  exceptional simplex or a Lawrence prism 
with $h_1 + \cdots + h_n \geq 2$ 
has degree $1$.
\label{d1}
\end{proposition}

\begin{proof}  If $\Delta$ is an $n$-dimensional Lawrence prism then, 
by \ref{Cay}, we have $\deg(\Delta) \leq 1$. On the other
hand, $\Vol(\Delta) = h_1 + \ldots + h_n$. Thus, we have
\[  (1-t)^{n+1} P(\Delta, t) = 1 + ( h_1 + \ldots + h_n -1)t \]
and  $\deg(\Delta) = 1$ for $h_1 + \cdots + h_n \geq 2$. 

If $\Delta_0$ is an exceptional triangle then, 
by  a direct computation, one obtains 
\[ (1-t)^3 P(\Delta_0, t) = 1 + 3t.  \]
For an arbitrary exceptional simplex $\Delta$, the statement follows 
from \ref{Cay}. 
\end{proof}

It turns out that the converse statement to \ref{d1} is true: 

\begin{theorem}[Main Theorem] 
Let $\Delta$ be lattice polytope of arbitrary dimension 
$n$. Then
$\deg(\Delta) \leq 1$ if and only if $\Delta$ is an exceptional simplex 
or a Lawrence prism.\label{maintheorem}
\end{theorem}

For $n=2$ this statement was proved independently 
by Arkinstall \cite{Ark80}, Khovanskii 
\cite[Sect.5]{Kho97}, Koelman \cite[Sect.4.1]{Koe91} 
and Schicho \cite[Thm.3.2]{Sch03}.\\

\smallskip

The starting point for the proof of this theorem 
is the following observation concerning the lattice points in $\Delta$:

\begin{lemma} 
Let $\Delta$ be a lattice polytope with  $\deg(\Delta) \leq 1$. 
Then the following statements are equivalent:
\begin{enumerate}
\item Lattice points in $\Delta$ are only vertices;
\item There is  no lattice point strictly between 
two lattice points in $\Delta$;
\item There is no lattice point strictly between 
two vertices of $\Delta$.
\end{enumerate}
\label{mainlemma}
\end{lemma}

\begin{proof} Obviously, one has 
$(1) \Rightarrow (2) \Rightarrow (3)$. By \ref{lat-poin},  
every lattice point in $\Delta$ lies on an edge. Thus, 
$(3) \Rightarrow (1)$.  
\end{proof}

\begin{definition}{\rm 
We call a lattice polytope $\Delta$ of degree $\leq 1$ 
 {\em narrow}, if it satisfies the   
equivalent statements in \ref{mainlemma}.
}
\end{definition}

The proof of the main theorem \ref{maintheorem} splits into two parts, 
depending on whether $\Delta$ is narrow or not.

\section{Classification of narrow lattice polytopes}

First we note that any lattice subpolytope of a narrow polytope is narrow. 

In the case of a simplex the situation is simple:

\begin{lemma}[Simplex-Lemma]
Let $S \subset \MR$ be an $n$-dimensional lattice simplex. Then $S$ is a basic simplex if and 
only if $S$ is narrow.
\label{basiclemma}
\end{lemma}

\begin{proof}

Let $S$ be narrow. By \ref{ele} we get $\Vol(\Delta) = 
\card{\Delta \cap M} - n = 1$.

\end{proof}

The next case to consider is that of a circuit, here very special 
lattice parallelograms naturally turn up. For this we fix some notation:

\begin{definition}{\rm Let $x_1,x_2,x_3,x_4$ be lattice points. 
If $\conv(x_1,x_2,x_3,x_4)$ is a two-dimensional $4$-gon $P$ 
with vertices $x_1,x_2,x_3,x_4$ that satisfy 
$x_1+x_4=x_2+x_3$ and $P$ contains no other lattice points except its vertices, 
then we call $x_1+x_4=x_2+x_3$ a {\em parallelogram relation} 
and $P$ a {\em narrow parallelogram}. In this case any 
subset of $\V(P)$ with three elements forms a 
two-dimensional affine lattice basis. The following 
figure illustrates this situation:

\centerline{\epsfxsize1in\epsffile{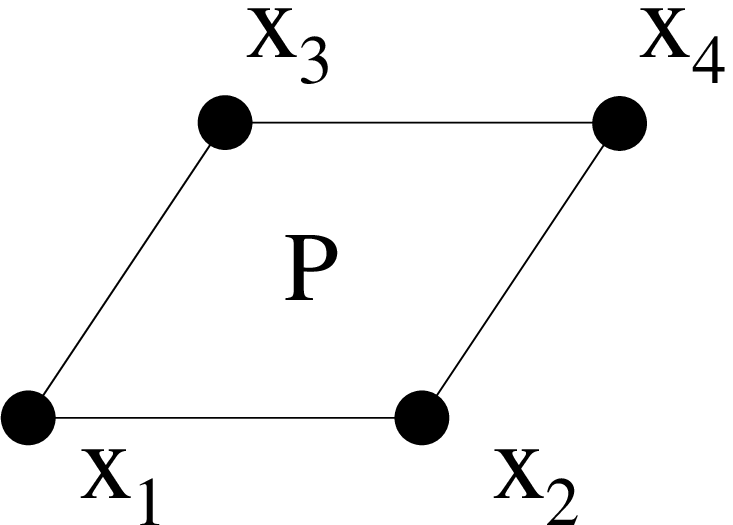}}

}
\end{definition}

\begin{lemma}[Circuit-Lemma]
Let $W \subset M$ be the set of $n+1$ vertices of a basic simplex $S$. 

Let $v \in M$, $v \not\in W$, such that 
$\conv(v,w) \cap M = \{v,w\}$ for any $w \in W$. 

If $\conv(v,W)$ has degree $\leq 1$, i.e., is narrow, then 
there is a parallelogram relation $v+w=w'+w''$ for some $w,w',w'' \in W$.
\label{circuitlemma}
\end{lemma}

\begin{proof}

Let $C$ be the unique circuit in $W \cup \{v\}$, that is, 
any proper subset of $C$ is affinely independent, 
however $C$ itself is not. Then $v \in C$, because $W$ is affinely 
independent. There is a unique affine relation
\[\sum_{i=1}^a c_i w_i = \sum_{j=1}^b d_j u_j,\]
with positive integers $c_i,d_j,a,b$, where $a+b = \card{C}$ and 
$C = \{w_1, \ldots, w_a, u_1,$\\$\ldots, u_b\}$. 

By \cite[Ch.7,Prop.1.2]{GKZ94} there are exactly two triangulations 
of $\Deltas := 
\conv(C)$ with vertices in $C$, namely 
$\{\conv(C \backslash \{w_i\}) \,:\, i = 1, \ldots, a\}$ and 
$\{\conv(C \backslash \{u_j\})$\\$\,:\, j = 1, \ldots, b\}$. 
By assumption and Simplex-Lemma \ref{basiclemma} any simplex in 
these triangulations has 
normalized volume one. Hence $\Vol(\Deltas) = a = b$. 
For $d := \dim(\Deltas) = \card{C} - 2 \geq 2$ Lemma 
\ref{ele} implies $$a + d = \Vol(\Deltas) + \dim(\Deltas) = 
\card{\Deltas \cap M} = \card{C} = 2 + d = a+b.$$ Thus  
$\Vol(\Deltas)=d=b=a=2$. This proves the statement.

\end {proof}

Let $e_0, \ldots, e_n$ be an affine lattice basis of $M$. 
Choose a lattice point $e'_l := e_l + e_n - e_0$ for some fixed 
$l \in \{1, \ldots, n-1\}$ so that 
$P_1 := \conv(e_l,e'_l,e_0,e_n)$ is a narrow parallelogram.

Choose another lattice point $v$ satisfying a parallelogram relation 
$v + e_c = e_a + e_b$ for some $a,b,c \in \{0, \ldots, n\}$, i.e., 
$P_2 := \conv(v,e_a,e_b,e_c)$ is another narrow parallelogram. 
Assume that 
$v - e_a \not= \pm (e_n - e_0)$ and $v - e_b \not= \pm (e_n - e_0)$.

\begin{lemma}[Parallelogram-Lemma] 
In the above situation, 
if the polytope $P := \conv(P_1,P_2)$ is narrow, then 
$P$ is a three-dimensional Lawrence prism, and we may assume 
(by possibly interchanging $a$ and $b$) that $P$ looks like

\centerline{\epsfxsize3in\epsffile{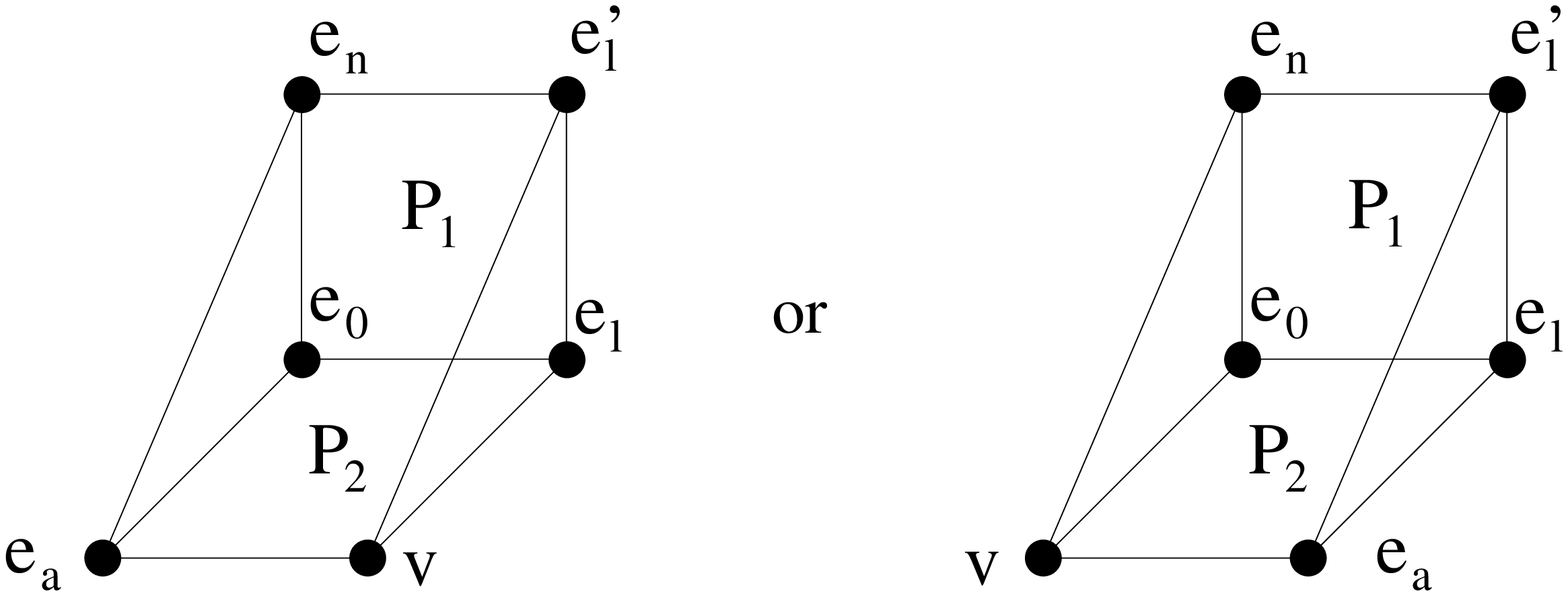}}

In both cases $e_0,e_l$ are the vertices of a common edge of the 
parallelograms $P_1$ and $P_2$.
\end{lemma}

\begin{proof}

If $\{a,b,c\} \cap \{0,l,n\} = \emptyset$, then 
$\dim(P) = 5$ and 
$$v + e_c + e_l + e_n = 4 \left(\frac{v + e_a + 
e_b + e_c + e_l + e'_l + e_0 + e_n}{8}\right)$$ 
is a lattice point in the interior of $4 P$, a contradiction.

If $\card{\{a,b,c\} \cap \{0,l,n\}} = 1$, then $\dim(P) = 4$. By symmetry 
we may assume $e_c = e_n$, so 
$$v + e_c + e_l = 3 \left(\frac{v + e_a + e_b + e_l + e'_l + e_0}{6}\right)$$
is a lattice point in the interior of $3 P$, a contradiction.

If $\card{\{a,b,c\} \cap \{0,l,n\}} = 2$, then $\dim(P) = 3$. 
Assume that $P_1$ and $P_2$ do not have a common edge. 
Then 
it is easy to see that the middle point $(e_l + e_n)/2$ of $P_1$ 
or the middle point $(v + e_c)/2$ of $P_2$ 
is in the interior of $P$, a contradiction.

If $\{a,b,c\} = \{0,l,n\}$, then $\dim(P) = 2$, and we immediately see 
that $P$ is not narrow, a contradiction.
\end{proof}

Now we can give the first part of the proof of Theorem \ref{maintheorem}.

\begin{proof}[Proof of Theorem \ref{maintheorem}, if $\Delta$ is narrow]

We may assume $\dim(\Delta) \geq 3$. The proof proceeds by induction 
on the number of lattice points in $\Delta$.

Due to the Simplex-Lemma \ref{basiclemma} we may assume that $\Delta$ 
is not a simplex, so 
there exists a vertex $v \in \V(\Delta)$ such that 
$\Deltas := \conv((\Delta \cap M)\backslash\{v\})$ is $n$-dimensional. 
By induction hypothesis $\Deltas$ is a 
Lawrence prism 
with respect to an affine lattice basis $e_0, \ldots, e_n$ and direction 
$e_n-e_0$. 

By the Circuit-Lemma \ref{circuitlemma} we get a parallelogram relation 
$v + e_c = e_a + e_b$ 
for $a,b,c \in \{0, \ldots, n\}$ pairwise different.

If $\card{\Deltas \cap M} = n+1$, then $\Delta$ is already a Lawrence 
prism, so we can assume that there is 
another vertex $e'_l \in \V(\Deltas)$, hence $e'_l = e_l + e_n - e_0$ for 
$l \in \{1, \ldots, n-1\}$. 

From now on we assume that $\Delta$ is not a Lawrence prism, so 
$v \not= e_a \pm (e_n - e_0)$ 
and $v \not= e_b \pm (e_n - e_0)$. 

In this situation the above Parallogram-Lemma yields that the convex 
hull $P$ of 
$e_l,e'_l,e_0,e_n,v,e_a,e_b,e_c$ is a three-dimensional Lawrence prism, 
and moreover we can assume that there are, 
as pictured, only two possibilities for $P$. In any case $e_0,e_l$ are 
vertices of the parallelogram $P_2 := \conv(v,e_a,e_b,e_c)$.

Since we assumed that 
$\Delta$ is not a Lawrence prism, there has to exist yet another vertex 
$e'_k := e_k + e_n - e_0\in \V(\Deltas)$ 
for $k \in \{1, \ldots, n-1\}\backslash\{l\}$. 

To finish the proof we apply the Parallelogram-Lemma again for $k$ 
instead of $l$, and get therefore 
that also $e_0,e_k$ are the vertices of an edge of the parallelogram $P_2$. 
Thus $Q := \conv(v,e_0,e_n,e_l,e'_l,e_k,e'_k)$ has to look like

\centerline{\epsfxsize1.25in\epsffile{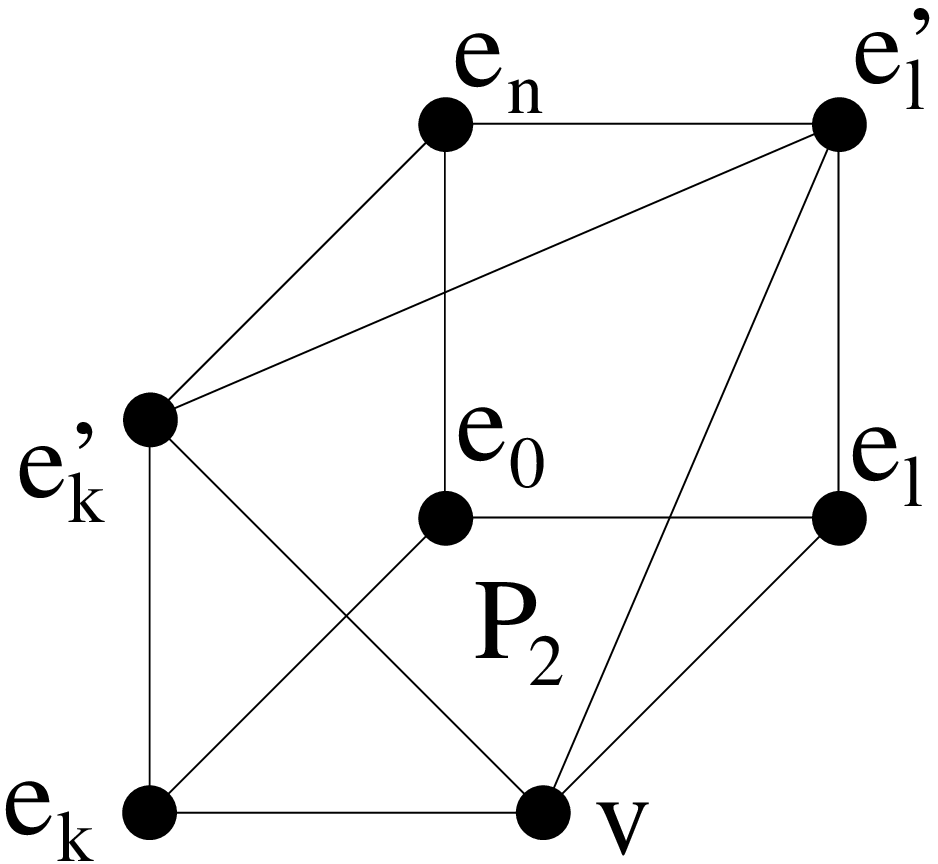}}

Obviously $(v + e_n)/2$ is a lattice point in the interior of $Q$, a 
contradiction.

\end{proof}

\section{Lattice polytopes of degree $1$ which are not narrow.}

\begin{definition}
{\rm Let  $\Delta \subseteq \MR$ be an $n$-dimensional lattice polytope.  
An edge of a lattice polytope is called {\em long}, if it contains more 
than two lattice points.}
\end{definition}

In particular for $\deg(\Delta) \leq 1$, Lemma \ref{mainlemma} implies 
that there exists a long edge if and only if $\Delta$ is not narrow.

\pagebreak
\begin{lemma}[Long-Edge-Lemma]
Let $\deg(\Delta) \leq 1$, and let $E$ be a long edge of $\Delta$.

\begin{enumerate}
\item If $E'$ is another long edge of $\Delta$, then $\conv(E,E')$ has 
dimension two, and is either an exceptional triangle or 
a Lawrence prism with parallel edges $E$ and $E'$.
\item If $D$ is an exceptional triangle contained in $\Delta$, then $E$ 
has to be an edge of $D$.
\item If $P$ is a narrow parallelogram contained in $\Delta$, and $Q := 
\conv(E,P)$ is not two-dimensional, then 
$\dim(Q)=3$ and $E \cap \V(P) = \emptyset$.
\end{enumerate}
\label{linelemma}
\end{lemma}

\begin{proof}

(1) If $u$ (resp. $u'$) is an interior lattice point of $E$ (resp. $E'$), 
then $(u+u')/2$ is a lattice point in the interior of $\conv(E,E')$.

(2) Applying (1) to all edges $E'$ of $D$ yields $\dim(\conv(E,D)) = 2$.

(3) Let $E = \conv(v,v')$ and $w,u \in E \cap M$ with $(v+w)/2 = u$. 
Let $x_1 + x_4 = x_2 + x_3$ be the parallelogram relation of $P$.

Since 
$$x_1 + x_4 + u = 3 \left(\frac{x_1 + x_2 + x_3 + x_4 + v + w}{6}\right)$$
is a lattice point in the interior of $3 Q$, 
we have $\dim(Q) = 3$. It is now easy to see that if $E \cap \V(P) \not= 
\emptyset$, then 
$(u + x_i)/2$ is a lattice point in the interior of $Q$ for some $i \in \{1, \ldots, 4\}$, a contradiction.

\end{proof}

This lemma is now applied to show that a long edge determines the 
direction of a Lawrence prism or otherwise lies in an exceptional triangle:

\begin{lemma}[Direction-Lemma]

Let $\deg(\Delta) \leq 1$. Let $e_0, \ldots, e_n \in \Delta$ be an affine 
lattice basis of $M$ 
such that $e_0$ and $e_1$ lie on a long edge $E$ of $\Delta$ 
having $e_0$ as a vertex and $e_1$ as an 
interior lattice point. Let $S := \conv(e_0, \ldots, e_n)$. The following 
figure illustrates the situation for $n=3$:

\centerline{\epsfxsize2in\epsffile{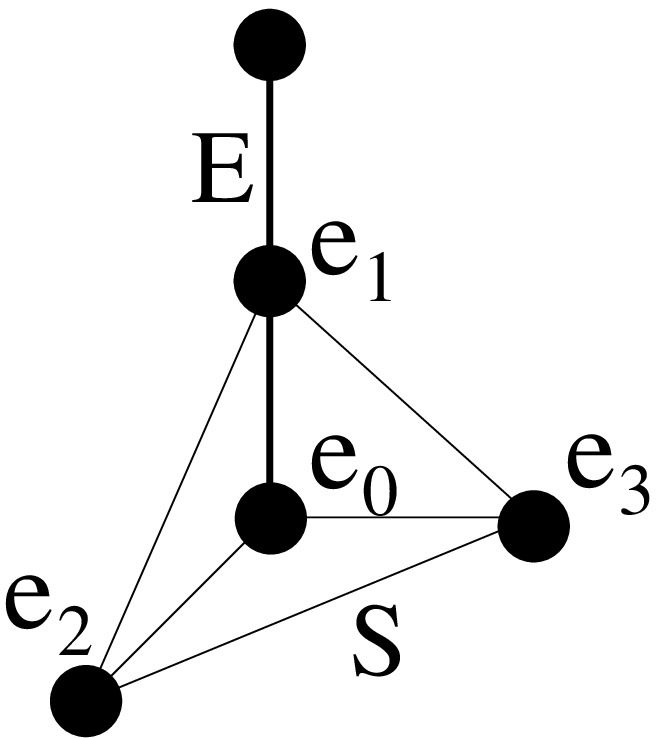}}

If $v$ is a vertex of $\Delta$, then $v \in \V(S) + \Z (e_1 - e_0)$ or 
$v$ is contained in an exceptional triangle in $\Delta$.
\end{lemma}

\begin{proof}

We may assume that $v$ is not contained in an exceptional triangle. 

If there exists some $w \in \V(S)$ so that $E' := \conv(v,w)$ is a long edge, then 
part (1) of the Long-Edge-Lemma implies that $E$ and $E'$ are parallel, hence 
$v \in \V(S) + \Z (e_1 - e_0)$. 

Otherwise $\conv(v,w) \cap M = \{v,w\}$ for any $w \in \V(S)$. We may now 
apply the Circuit-Lemma to $\V(S)$ and $v$ to get a narrow 
parallelogram $P$ with vertex $v$ and the other three vertices in $\V(S)$.
Let $Q := \conv(E,P)$. Assume $\dim(Q) > 2$. Then part (3) of the 
Long-Edge-Lemma implies $\dim(Q) = 3$ and $\V(P) \subseteq \{v,e_2, \ldots, 
e_n\}$. Hence 
$\V(P)\backslash\{v\}, e_0, e_1$ is affinely independent, so $\dim(Q) = 4$, 
a contradiction. 
Hence $\dim(Q) = 2$, thus $e_0,e_1 \in \V(P)$. Let $e_j \in \V(P)$ for $j 
\in \{2, \ldots, n\}$; we may assume $j=2$. 
If $v \not= e_2 \pm (e_1 - e_0)$, then $v = e_0 + e_1 - e_2$, so looking at 
the figure yields $e_1$ as an interior lattice point of the two-dimensional polygon $\conv(v,e,e_2)$, a contradiction.

\end{proof}

Now we can finish the proof of Theorem \ref{maintheorem}.

\begin{proof}[Proof of Theorem \ref{maintheorem}, if $\Delta$ is not narrow]

Let $\deg(\Delta) \leq 1$. We may assume $\dim(\Delta) \geq 3$.

Since $\Delta$ is assumed to be not narrow, we can find a long edge $E$ of 
$\Delta$ with vertices $e_0$ and $e'_1$, 
a lattice point $e_1$ in the interior of $E$ such that $e_1 - e_0$ is a 
primitive lattice point, and we find 
$e'_2, e_3, \ldots, e_n \in \V(\Delta)$ such that $e_0,e_1,e'_2, \ldots, e_n$ 
are the vertices of an $n$-dimensional simplex $S'$. 
So we are in the situation of the Direction-Lemma, however $e_0,e_1,e'_2,e_3, 
\ldots, e_n$ is a priori just an affine $\R$-basis of $\MR$.

If $\Delta$ contains no exceptional triangles, then by part (1) of the 
Long-Edge-Lemma 
there cannot be a long edge in $S'$, thus $S'$ is basic by the Simplex-Lemma 
\ref{basiclemma}, and the Direction-Lemma yields 
that $\Delta$ is a Lawrence prism with direction $e_1 - e_0$.

Hence we may assume that $D := \conv(e_0,e'_1,e'_2)$ is an exceptional triangle.

Let $e_2 := (e_0 + e'_2)/2 \in M$. By part (2) of the Long-Edge-Lemma {\em any} 
long edge of $\Delta$ is contained in $D$.
In particular the simplex $S := 
\conv(e_0,e_1,e_2, \ldots, e_n)$ is narrow, hence basic. 
This implies that the simplex $P := \conv(e_0,e'_1,e'_2,e_3, \ldots, e_n)$ is exceptional. 
$P$ is illustrated in the following figure (for $n=3$):

\centerline{\epsfxsize2in\epsffile{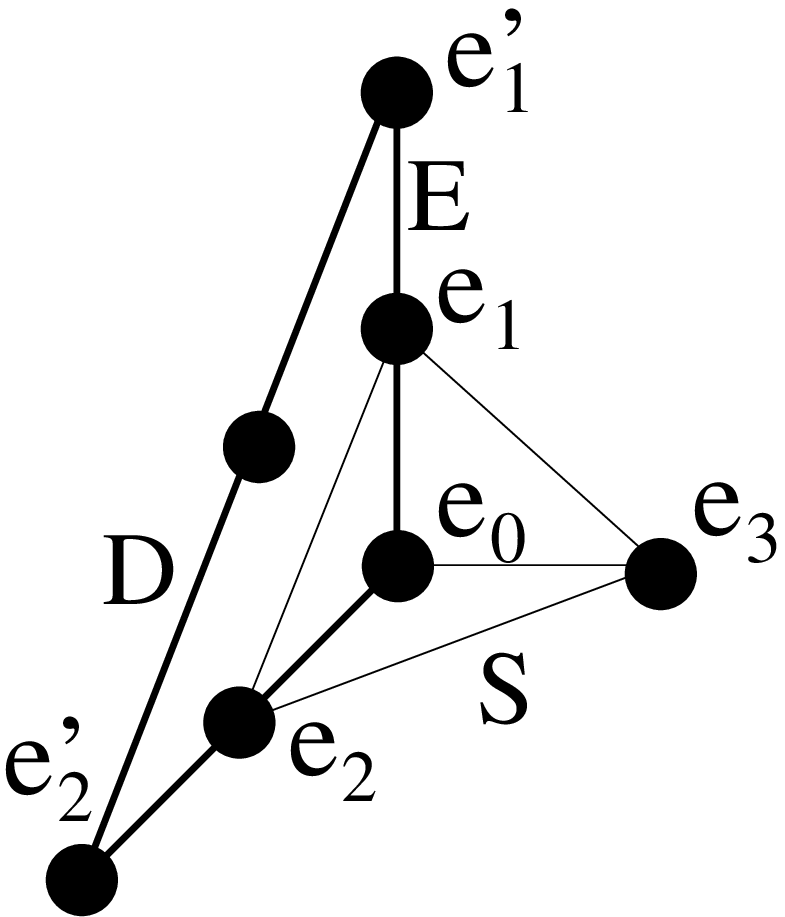}}

Assume there exists $v \in \V(\Delta)$ with $v \not\in P$. 

In particular $v$ cannot be the vertex of a long edge. Hence we can apply the Direction-Lemma to the long edge $P$ and the vertices of $S$ 
to get $v = e_j \pm (e_1 - e_0)$ for some 
$j \in \{2, \ldots, n\}$. If $j=2$, then either $v \in P$ or $e_2$ is an 
interior lattice point of the polygon $\conv(v,e'_2,e'_1,e_0)$, a 
contradiction. 
Hence $j > 2$. However applying the Direction-Lemma again, this time to 
the long edge $\conv(e'_2,e_0)$, yields 
$v \in  \V(S) + \Z (e_2 - e_0)$, a contradiction.

Thus $\Delta = P$.

\end{proof}

\section{Coherent triangulations and $A$-determinant}

As some immediate applications we determine for all lattice polytopes of degree 
$\leq 1$ their triangulations, secondary polytopes and 
principal $A$-determinants.

\begin{proposition}
Let $\deg(\Delta) \leq 1$. Then any lattice triangulation $T$ of $\Delta$ is 
{\em coherent}, i.e., there exists a concave $T$-piecewise linear function whose domains of 
linearity are precisely the (maximal) simplices of $T$. 

If $\Delta$ is exceptional, then there are $14$ lattice triangulations.

If $\Delta$ is the Lawrence prism with heights $h_1, \ldots, h_n$, then 
the lattice triangulations of $\Delta$ are in one-to-one correspondendence
 with words over the alphabet
\[\{x_{1,1},\ldots, x_{1,h_1}, \ldots, x_{n,1}, \ldots, x_{n,h_n}\}\]
where for each $i$ the sum of the second index of the letters in $\{x_{i,1}, 
\ldots, x_{i,h_i}\}$ has to be $h_i$. 
The number of lattice triangulations of $\Delta$ is
\[\sum\limits_{(l_1, \ldots, l_n) \in \{1, \ldots, h_1\} \times \cdots \times \{1, 
\ldots, h_n\}} 
\frac{(l_1 + \cdots + l_n)!}{{l_1}! \cdots {l_n}!} \prod\limits_{i=1}^n 
\binom{h_i-1}{l_i-1}.\]
\end{proposition}

The proof of this corollary is left to the reader 
(e.g., compare with \cite[Sect.3]{DHZ01} and \cite[Ch.7,Prop.3.4]{GKZ94}).

In the case of the classical discriminant we have $n=1$ and $\Delta = [0,h]$, 
and the above formula gives a well-know result: 
there exist exactly $2^{h-1}$ coherent triangulations of $\Delta$.

In the case of an $n$-dimensional Lawrence prism with heights $1, \ldots, 1$,
 that is, just the product of 
an $n-1$-dimensional basic simplex and a unit segment, the secondary polytope 
is the well-known 
{\em permutahedron}, cf. \cite[Ch.7,Sect.3]{GKZ94}. It is an $n-1$-dimensional polytope $P$ with $n!$ vertices 
that is {\em simple}, i.e., any vertex of $P$ is contained in only $\dim(P) = n-1$ edges. 
Here we prove the following generalization:

\begin{proposition}
Let $\deg(\Delta) \leq 1$. Then the secondary polytope $\Secf$ is a 
{\em simple} lattice polytope.

If $\Delta$ is exceptional, then $\Secf$ is the three-dimensional {\em associahedron}, cf. \cite[Ch.7,Sect.3B]{GKZ94}. 
The following figure is a \texttt{polymake}-visualization, cf. \cite{GJ05}:

\centerline{\epsfxsize1.25in\epsffile{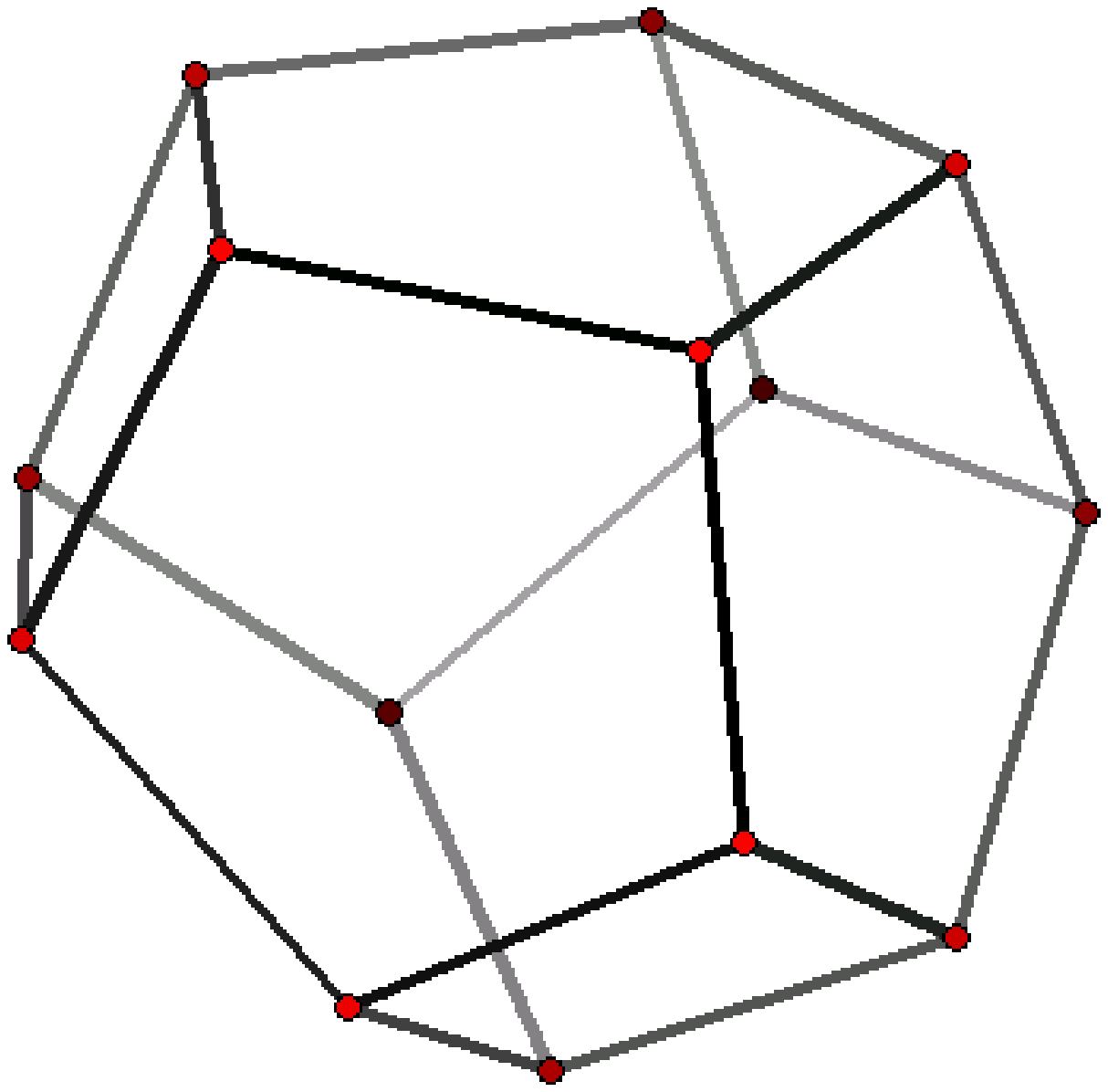}}

If $\Delta$ is the Lawrence prism 
with heights $h_1, \ldots, h_n$, then $\Secf$ is $(h_1 + \cdots + h_n - 1)$-dimensional.
\end{proposition}

\begin{proof}

For $\Delta$ exceptional it is enough to compute the secondary polytope for $n=2$, cf. \cite[Ch.13,Sect.1,E. Example]{GKZ94}.

So let $\Delta$ be a Lawrence prism. 
We get $\dim(\Secf) = \card{\Delta \cap M}  - n - 1 = (h_1 + 1) + \cdots (h_n + 1) - n - 1 
= h_1 + \cdots + h_n - 1$. 

By \cite[Ch.7,Thm.2.10]{GKZ94} two vertices of $\Secf$ are joined by an edge if 
and only if they are obtained by a modification along a circuit. So let $T$ be a triangulation with 
maximal simplices $T_1, \ldots, T_s$, sorted in ascending order with respect to their intersection with 
the affine line containing $(e_0 + \cdots + e_{n-1})/n$ and with direction $e_n - e_0$. Let $T'$ be another triangulation 
such that the associated vertices in $\Secf$ are joined by an edge. 

If $s=1$, then $\Delta = T_1$ is a simplex. Hence we may assume $\Delta = (e_0, \ldots, e_{n-1}, e_0 + h_1 (e_n - e_0))$. 
Therefore any lattice triangulation is induced by a lattice triangulation of $[0,h_1]$. 
However $\Sec([0,h_1])$ is combinatorially a cube by \cite[Ch.7,Prop.3.1]{GKZ94}, in particular simple.

So let $s > 1$, and we may assume 
$T_1 = \conv(e_0, e_1, \ldots, e_{n-1}, e_0 + k (e_n - e_0))$.

If $T_1 \cap T_2$ is also a simplex of $T'$, i.e., it stays unmodified under the modification of $T$, then $T'$ is either 
given by a lattice triangulation of the Lawrence prism 
$\conv(T_2,T_s)$ or by a lattice triangulation of $T_1$. By induction hypothesis there are 
$h_1 + \cdots + h_n - 2 = (h_1 + \cdots + h_n - k - 1) + (k - 1)$ many choices for such $T'$. 

Otherwise $T'$ is the modification of $T$ along the unique circuit contained in 
$\conv(T_1,T_2)$. This proves that $\Secf$ is simple.

\end{proof}

For $A := \Delta \cap $M the principal $A$-determinant was computed in \cite{GKZ94} in the cases that 
$\Delta$ is an exceptional triangle \cite[Ch.13,Sect.1,E. Example]{GKZ94} or a product 
of two basic simplices \cite[Ch.10,Example 1.3(b)]{GKZ94}. Here we note the formula 
for a Lawrence prism:

\begin{proposition}
Let $E_A(f)$ be the principal $A$-determinant for $A := \Delta \cap M$ and the polynomial 
$f = \sum_{\omega \in \Delta \cap M} a_\omega X^\omega$.

If $\Delta$ is the Lawrence prism with heights $h_0, \ldots, h_{n-1}$, then we can assume $e_0 = 0$ and 
$E_A(f) = f_0(x) + \sum_{i=1}^{n-1} y_i f_i(x)$ for $f_i(x) = \sum_{j=0}^{h_i} a_{i,j} x^j$. 

Then
\[E_A(f) = (\prod_{i=0}^{n-1} a_{i,0} a_{i,h_i}) \; (\prod_{i=0}^{n-1} \Discr_{f_i}) \; 
(\prod_{i < j} \Res(f_i,f_j)).\]
\end{proposition}

\begin{proof}

The degree of $E_A(f)$ is $(n+1) \Vol(\Delta)$, that is the degree of the polynomial on the right side. 
However by prime factorization \cite[Ch.10,Thm.1.2]{GKZ94} and the Cayley-Trick \cite[Ch.9,Prop.1.7]{GKZ94} 
the right side divides the left.

\end{proof}

The case where all heights equal one, i.e., $\Delta$ is a product of a basic simplex and a unit interval, 
was already explicitly treated in \cite[Ch.10,Example 1.5(b)]{GKZ94}. Here we illustrate this formula for 
$n=3$ with $h_0=h_1=1$, $h_2=2$ and coefficients enumerated as in the following figure:

\centerline{\epsfxsize1.25in\epsffile{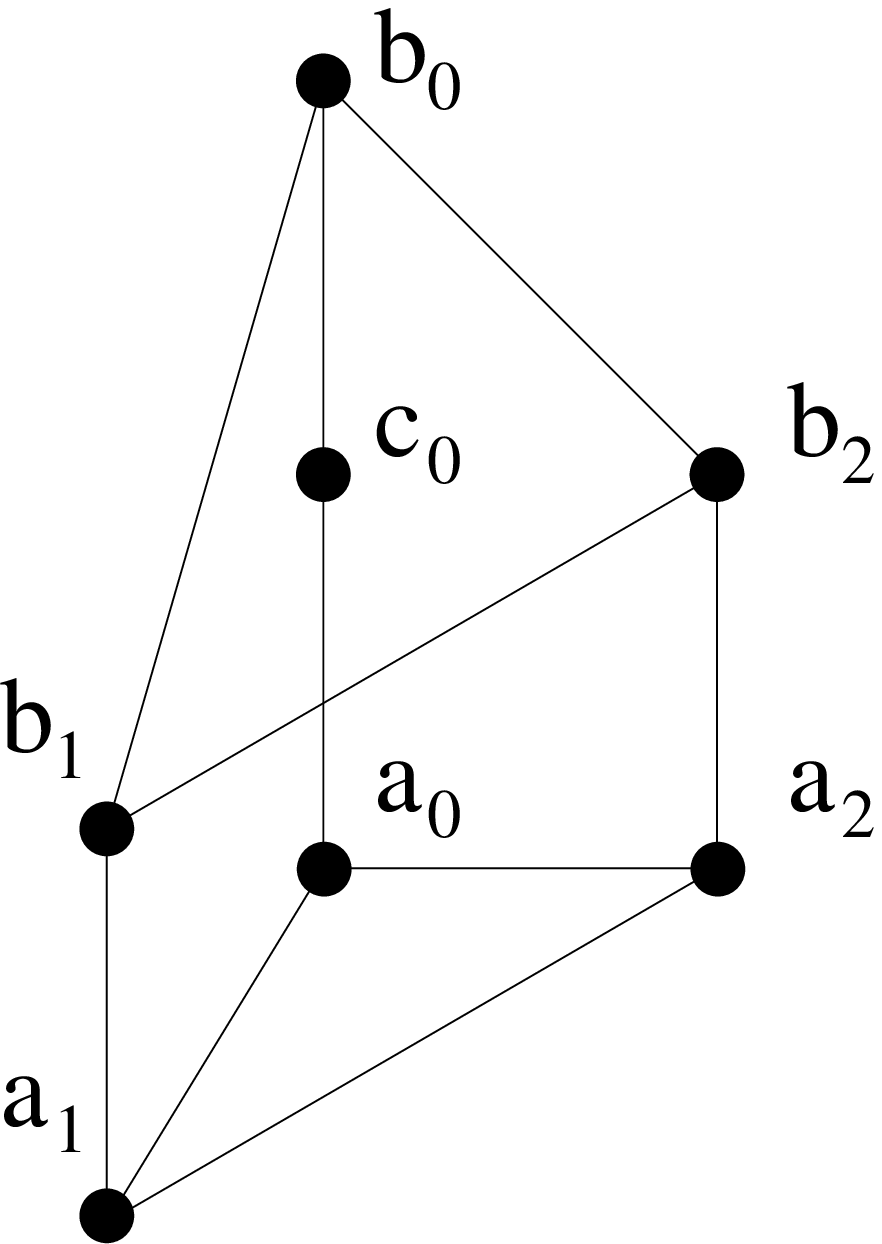}}

\vspace{-1ex}
Then
\begin{eqnarray*}E_A(f) &=&a_0 a_1 a_2 b_0 b_1 b_2 \,(4 a_0 b_0 - c_0^2)\\
&&(a_1 b_2 - a_2 b_1) (a_1^2 b_0 + b_1^2 a_0 - a_1 b_1 c_0) (a_2^2 b_0 + b_2^2 a_0 - a_2 b_2 c_0).
\end{eqnarray*}

\end{document}